\documentclass{commat}

\title[On the degree of approximation of continuous functions]{On the degree of approximation of continuous functions by a linear transformation of their Fourier series}

\author[Xh. Z. Krasniqi]{Xhevat Z. Krasniqi}
\affiliation{\address{Faculty of Education, University of Prishtina ``Hasan Prishtina'', Avenue ``Mother Theresa'' 5, 10000 Prishtina, Republic of Kosovo}
\email{xhevat.krasniqi@uni-pr.edu}}



\abstract{
In this paper, we have proved four theorems on the degree of approximation of continuous functions by matrix means of their Fourier series which is expressed in terms of the modulus of continuity and a non-negative mediate function.}

\firstpage{37}

\VOLUME{30}

\DOI{https://doi.org/10.46298/cm.9273}

\begin{paper}

\section{Introduction}

By $S_{n}(f;x)$ we denote the
$n$-th partial sum of the Fourier series of a $2\pi$-periodic continuous function $f(x)$ at $x$ and $\omega (\delta) = \omega (\delta, f)$ the modulus of continuity of $f$.

In short, we denote by $A: = (a_{n,k})$, $(k, n = 0, 1,\dots)$ a lower triangular infinite matrix of real numbers and the $A$-transform of $\{S_{n}(f;x)\}$ by
\[
T_{n,A}(f; x) : = \sum_{k = 0}^{n}a_{n,k}S_{k}(f; x),\quad (n = 0, 1, \dots).
\]

Throughout this paper we write $u = \mathcal{O}(v)$ if there exists a positive constant $C$ such that $u\leq Cv$, and $\|\cdot\|$ denotes the sup-norm.

P. Chandra~\cite{CH1}, \cite{CH2} was the first who estimated the deviation $T_{n,A}(f)-f$, in sup-norm, when
$\{a_{n,k}\}$ is a monotonic sequences with respect to $k$. Many years later, L. Leindler~\cite{L} replaced the monotonicity condition with the so-called \emph{rest bounded variation sequences} and \emph{``head bounded variation sequences''} (these will be recalled later in this paper).

Chandra's theorems are the following.

\begin{theorem}\label{theorem1}
Let $\{a_{n,k}\}$ satisfy the following conditions:
\begin{equation}\label{eq11}a_{n,k} \geq 0\qquad \mbox{and}\qquad \sum_{k = 0}^{n}a_{n,k} = 1,
\end{equation}
\begin{equation}\label{eq12}a_{n,k} \leq a_{n,k+1}\quad (k = 0,1,\dots,n-1; n = 0,1,\dots).
\end{equation}
Suppose $\omega (t)$ is such that
\begin{equation}\label{eq13}\int _{u}^{\pi}t^{-2}\omega (t)dt = \mathcal{O}\left(H(u)\right)\quad(u\to +0),
\end{equation}
where $H(u)\geq 0$ and
\begin{equation}\label{eq14}\int _{0}^{t}H(u)du = \mathcal{O}\left(tH(t)\right)\quad(t\to +0).
\end{equation}
Then
\begin{equation}\label{eq15}\nonumber \|T_{n,A}(f)-f\| = \mathcal{O}\left(a_{n,n}H(a_{n,n})\right).
\end{equation}
\end{theorem}

\begin{theorem}\label{theorem2}
Let~\eqref{eq11},~\eqref{eq12} and~\eqref{eq13} hold. Then
\begin{equation}\label{eq16}\nonumber \|T_{n,A}(f)-f\| = \mathcal{O}\left(\omega (\pi / n)\right)+\mathcal{O}\left(a_{n,n}H(\pi / n)\right).
\end{equation}

If, in addition, $\omega (t)$ satisfies~\eqref{eq14} then
\begin{equation}\label{eq17}\nonumber \|T_{n,A}(f)-f\| = \mathcal{O}\left(a_{n,n}H(\pi / n)\right).
\end{equation}
\end{theorem}

\begin{theorem}\label{theorem3}
Let us assume that~\eqref{eq11} and
\begin{equation}\label{eq18}a_{n,k} \geq a_{n,k+1}\quad (k = 0,1,\dots,n-1; n = 0,1,\dots)
\end{equation}
hold.
Then
\begin{equation}\label{eq19}\nonumber \|T_{n,A}(f)-f\| = \mathcal{O}\left(\omega (\pi / n)+\sum_{k = 1}^{n}k^{-1}\omega (\pi / k) \sum_{r = 0}^{k+1} a_{n,r})\right).
\end{equation}
\end{theorem}

\begin{theorem}\label{theorem4}
Let~\eqref{eq11},~\eqref{eq13},~\eqref{eq14} and~\eqref{eq18} hold. Then
\begin{equation}\label{eq20}\nonumber \|T_{n,A}(f)-f\| = \mathcal{O}\left(a_{n,0}H(a_{n,0})\right).
\end{equation}
\end{theorem}
\

To obtain his results, L. Leindler~\cite{L} used the following classes of sequences.

A sequence $\textbf{c} : = \{c_{n}\}$ of nonnegative numbers tending to zero is called of Rest Bounded Variation, or briefly $\textbf{c} \in RBV S$, if it has the property
\begin{equation}\label{eq21}\nonumber \sum_{n = m}^{\infty}|c_{n}-c_{n+1}|\leq K(\textbf{c})c_{m}
\end{equation}
for all natural numbers $m$, where $K(\textbf{c})$ is a constant depending only on $\textbf{c} $.

A sequence $\textbf{c} : = \{c_{n}\}$ of nonnegative numbers will be called of Head Bounded Variation, or briefly $\textbf{c}\in HBV S$, if it has the property
\begin{equation}\label{eq22}\nonumber \sum_{n = 0}^{m-1}|c_{n}-c_{n+1}|\leq K(\textbf{c})c_{m}
\end{equation}
for all natural numbers $m$, or only for all $m\leq N$ if the sequence $\textbf{c} $ has only finite nonzero terms, and the last nonzero term is $c_{N}$.

To go further it is clear that condition $0<K(\textbf{c})\leq K$ needs to be assumed, where $K$ is a positive constant.

Assuming that for all $n$ and $0\leq m\leq n$
\begin{equation}\label{eq23}\sum_{k = m}^{\infty}|a_{n,k}-a_{n,k+1}|\leq Ka_{n,m}
\end{equation}
and
\begin{equation}\label{eq24}\sum_{k = 0}^{m-1}|a_{n,k}-a_{n,k+1}|\leq Ka_{n,m}
\end{equation}
hold, where $K$ is a positive constant, Leindler~\cite{L} proved the following:

\begin{theorem}\label{theorem5}
The statements of Theorems~\ref{theorem1},~\ref{theorem2},~\ref{theorem3} and~\ref{theorem4} hold with~\eqref{eq24}
in place of~\eqref{eq12}, and with~\eqref{eq23} in place of~\eqref{eq18}, respectively; naturally maintaining all the other assumptions.
\end{theorem}

These results are extended further in~\cite{XhK}, then they are generalized in~\cite{XhK1}, and are treated again in~\cite{XhK2} using some other means (see also~\cite{BD}). Very recently, W.~Lenski and B.~Szal~\cite{LS}, used (for trigonometric approximation in the variable space $L^{p(x)}$ and $m = n$) the following condition
\begin{equation}\label{star}
\sum_{k = 0}^{m-1}{(k+1)}^{\beta}\left|\frac{a_{n,k}}{{(k+1)}^{\beta}}-\frac{a_{n,k+1}}{{(k+2)}^{\beta}}\right| = \mathcal{O}(a_{n,m}),\qquad(\beta \geq 0),
\end{equation}

assumed on the sequence $\{a_{n,k}\}$, $0\leq m\leq n$. To my best knowledge, such a condition appears for the first time in the mathematical literature, which for $\beta = 0$ coincides with condition~\eqref{eq24}.

Motivated by the above condition (for $a_{n,m} = \frac{1}{n+1}$ and $m = n$ this condition has been introduced in~\cite{LS}) we introduce another new condition
\begin{equation}\label{starstar}
\sum_{k = m}^{\infty}{(k+1)}^{\beta}\left|\frac{a_{n,k}}{{(k+1)}^{\beta}}-\frac{a_{n,k+1}}{{(k+2)}^{\beta}}\right| = \mathcal{O}(a_{n,m}),\qquad (\beta \geq 0),
\end{equation}

on the sequence $\{a_{n,k}\}$, $0\leq m\leq n$, which for $\beta = 0$ coincide with condition~\eqref{eq23}.

\begin{remark}
Note that for $\beta = 0$ these two conditions coincides with~\eqref{eq24} and~\eqref{eq23}, respectively.
Also, it was shown in~\cite{LS} that if $\{{(k+1)}^{-\beta}a_{n,k}\}\in HBVS$, then condition $\eqref{star}$ holds true. Moreover, it easy to verify that condition $\eqref{starstar}$ implies $\{{(k+1)}^{-\beta}a_{n,k}\}\in RBVS$. Furthermore, some examples, which supports these conditions, are given there as well.
\end{remark}

Now we are concerned weather these two conditions on the sequence $\{a_{n,k}\}$ can be used in Theorems~\ref{theorem1}--\ref{theorem4}, instead of~\eqref{eq12} and~\eqref{eq15}, which is in fact the aim of this paper.

\vspace{5mm}

\section{Auxiliary Results}

The next lemmas play a helpful role for the proof of our results.

\begin{lemma}
\cite{CH2}\label{theorem7}
If~\eqref{eq13} and~\eqref{eq14} hold then
\begin{equation}\label{eq26}\nonumber \int_{0}^{v}t^{-1}\omega(t)dt = \mathcal{O}\left(vH(v)\right)\quad (v\to +0).
\end{equation}
\end{lemma}

\begin{lemma}\label{theorem6}
If~\eqref{eq13} and~\eqref{eq14} hold then
\begin{equation}\label{eq25}\nonumber \int_{0}^{\pi/m}\omega(t)dt = \mathcal{O}\left(m^{-2}H(\pi/m)\right).
\end{equation}
\end{lemma}

\begin{proof}
Using Lemma~\ref{theorem7} we have
\begin{align*}
\int_{0}^{\pi/m}\omega(t)dt &= \int_{0}^{\pi/m}\frac{t\omega(t)}{t}dt\leq \frac{\pi}{m}\int_{0}^{\pi/m}\frac{\omega(t)}{t}dt\\
&= \frac{\pi}{m}\mathcal{O}\left(\frac{\pi}{m}H\left(\frac{\pi}{m}\right)\right) = \mathcal{O}\left(m^{-2}H(\pi/m)\right),
\end{align*}
which completes the proof.
\end{proof}

\begin{lemma}
\cite{LS}\label{theorem8}
If $\beta \geq 0$ and $0<t\leq \pi$, then
\[
\left|\sum_{j = 0}^{m}{(j+1)}^{\beta}\frac{\sin\left(j+\frac{1}{2}\right)t}{2\sin\frac{t}{2}}\right|\leq \frac{\pi^2 {(m+1)}^{\beta}}{t^2}.
\]
\end{lemma}

\begin{lemma}\label{theorem9} Let $\beta \geq 0$ be a real number.
If for a fixed $n$, the sequence $\{a_{n,k}\}$ satisfies the condition
\[
\sum_{k = m}^{\infty}{(k+1)}^{\beta}\left|\frac{a_{n,k}}{{(k+1)}^{\beta}}-\frac{a_{n,k+1}}{{(k+2)}^{\beta}}\right| = \mathcal{O}(a_{n,m}),
\]
for $0\leq m\leq n$, then uniformly in $0 < t\leq \pi$,
\begin{equation}\label{eq27}\nonumber
|K_n(t)|: = \left|\sum_{k = 0}^{n}a_{n,k}\frac{\sin\left(k+\frac{1}{2}\right)t}{2\sin \frac{t}{2}}\right| = \mathcal{O}\left(t^{-1}A_{n,\tau}\right),
\end{equation}
where $\displaystyle A_{n,\tau}: = \sum_{r = 0}^{\tau}a_{n,r}$ and $\tau$ denotes the integer part of $\,\frac{\pi}{t}$.

If $\{a_{n,k}\}$ satisfies the condition
\[
\sum_{k = 0}^{m-1}{(k+1)}^{\beta}\left|\frac{a_{n,k}}{{(k+1)}^{\beta}}-\frac{a_{n,k+1}}{{(k+2)}^{\beta}}\right| = \mathcal{O}(a_{n,m}),
\]
for $0\leq m\leq n$, then
\begin{equation}\label{eq28}
\nonumber |K_n(t)| = \mathcal{O}\left(\frac{a_{n,n}}{t^2}\right).
\end{equation}
\end{lemma}

\begin{proof}
Assume that $n\geq \tau$. Applying summation by parts, the inequalities $|\sin \gamma|\leq 1$ and $\pi \sin \gamma\geq 2\gamma$, ($\gamma \in[0,\pi/2]$), we have
\begin{align*}
\left|K_n(t)\right|&\leq \sum_{k = 0}^{\tau}\frac{a_{n,k}}{2\sin \frac{t}{2}}
+\left|\sum_{k = \tau+1}^{n}a_{n,k}\frac{\sin\left(k+\frac{1}{2}\right)t}{2\sin \frac{t}{2}}\right|\\
&\leq \frac{\pi}{2t}\sum_{k = 0}^{\tau}a_{n,k}
+\left|\sum_{k = \tau+1}^{n}\frac{a_{n,k}}{{(k+1)}^{\beta}}{(k+1)}^{\beta}\frac{\sin\left(k+\frac{1}{2}\right)t}{2\sin \frac{t}{2}}\right|\\
&\leq \frac{a_{n,\tau+1}}{{(\tau+2)}^{\beta}}\left|\sum_{j = 0}^{\tau}{(j+1)}^{\beta}\frac{\sin\left(j+\frac{1}{2}\right)t}{2\sin \frac{t}{2}}\right|\\
&\quad+\frac{\pi}{2t}\sum_{k = 0}^{\tau}a_{n,k}+\frac{a_{n,n}}{{(n+1)}^{\beta}}\left|\sum_{j = 0}^{n}{(j+1)}^{\beta}\frac{\sin\left(j+\frac{1}{2}\right)t}{2\sin \frac{t}{2}}\right|\\
&\quad+\left|\sum_{k = \tau+1}^{n-1}\left(\frac{a_{n,k}}{{(k+1)}^{\beta}}-\frac{a_{n,k+1}}{{(k+2)}^{\beta}}\right)\sum_{j = 0}^{k}{(j+1)}^{\beta}\frac{\sin\left(j+\frac{1}{2}\right)t}{2\sin \frac{t}{2}}\right|.
\end{align*}

Now, using Lemma~\ref{theorem8} and our assumption, we get
\begin{align*}
\left|K_n(t)\right|
&\leq \frac{a_{n,\tau+1}}{{(\tau+2)}^{\beta}}\frac{\pi^2 {(\tau+1)}^{\beta}}{t^2}+\frac{\pi}{2t}\sum_{k = 0}^{\tau}a_{n,k}\\
&\quad+\frac{a_{n,n}}{{(n+1)}^{\beta}}\frac{\pi^2 {(n+1)}^{\beta}}{t^2}
+\sum_{k = n}^{\infty}\left|\frac{a_{n,k}}{{(k+1)}^{\beta}}-\frac{a_{n,k+1}}{{(k+2)}^{\beta}}\right|\frac{\pi^2 {(k+1)}^{\beta}}{t^2}\\
&= \frac{\pi}{2t}\sum_{k = 0}^{\tau}a_{n,k}+\frac{\pi^2}{t^2}\mathcal{O}(a_{n,\tau}+a_{n,n}).
\end{align*}

But (by assumption) we also get
\begin{align*}
a_{n,n} &= \sum_{k = n}^{\infty}{(k+1)}^{\beta}\left|\frac{a_{n,k}}{{(k+1)}^{\beta}}-\frac{a_{n,k+1}}{{(k+2)}^{\beta}}\right|\\
&\leq \sum_{k = \tau}^{\infty}{(k+1)}^{\beta}\left|\frac{a_{n,k}}{{(k+1)}^{\beta}}-\frac{a_{n,k+1}}{{(k+2)}^{\beta}}\right|
= \mathcal{O}(a_{n,\tau}).
\end{align*}

Therefore,
\begin{equation*}
\left|K_n(t)\right|
= \frac{\pi}{2t}\sum_{k = 0}^{\tau}a_{n,k}+\frac{\pi}{t}\mathcal{O}((\tau+1)a_{n,\tau}) = \mathcal{O}\left(t^{-1}\sum_{k = 0}^{\tau}a_{n,k}\right),
\end{equation*}
since for $0\leq k\leq \tau$ we have
\begin{align*}
a_{n,\tau}&\leq \sum_{i = \tau}^{\infty}{(i+1)}^{\beta}\left|\frac{a_{n,i}}{{(i+1)}^{\beta}}-\frac{a_{n,i+1}}{{(i+2)}^{\beta}}\right|\\
&\leq \sum_{i = k}^{\infty}{(i+1)}^{\beta}\left|\frac{a_{n,i}}{{(i+1)}^{\beta}}-\frac{a_{n,i+1}}{{(i+2)}^{\beta}}\right|
= \mathcal{O}(a_{n,k}).
\end{align*}

Similarly, we have obtained
\begin{align*}
\left|K_n(t)\right|&\leq \frac{a_{n,n}}{{(n+1)}^{\beta}}\left|\sum_{j = 0}^{n}{(j+1)}^{\beta}\frac{\sin\left(j+\frac{1}{2}\right)t}{2\sin \frac{t}{2}}\right|\\
&\quad+\left|\sum_{k = 0}^{n-1}\left(\frac{a_{n,k}}{{(k+1)}^{\beta}}-\frac{a_{n,k+1}}{{(k+2)}^{\beta}}\right)\sum_{j = 0}^{k}{(j+1)}^{\beta}\frac{\sin\left(j+\frac{1}{2}\right)t}{2\sin \frac{t}{2}}\right|\\
&\leq \frac{a_{n,n}}{{(n+1)}^{\beta}}\frac{\pi^2 {(n+1)}^{\beta}}{t^2}
+\sum_{k = 0}^{n-1}\left|\frac{a_{n,k}}{{(k+1)}^{\beta}}-\frac{a_{n,k+1}}{{(k+2)}^{\beta}}\right|\frac{\pi^2 {(k+1)}^{\beta}}{t^2}\\
&= \mathcal{O}\left(\frac{a_{n,n}}{t^2}\right).
\end{align*}

The proof is now complete.
\end{proof}

\section{Main Results}
We start with the following.

\begin{theorem}\label{theorem10}
Let $\{a_{n,k}\}$ satisfy the conditions~\eqref{eq11} and
\[
\sum_{k = 0}^{m-1}{(k+1)}^{\beta}\left|\frac{a_{n,k}}{{(k+1)}^{\beta}}-\frac{a_{n,k+1}}{{(k+2)}^{\beta}}\right| = \mathcal{O}(a_{n,m}),
\]
for all $m = 0,1,\dots,n;\, (n = 0,1,\dots)$.
Suppose $\omega (t)$ is such that~\eqref{eq13} and~\eqref{eq14} hold, then
\begin{equation}\label{eq55}
\|T_{n,A}(f)-f\| = \mathcal{O}\left(a_{n,n}H\left(a_{n,n}\right)\right).
\end{equation}
\end{theorem}
\begin{proof}
Writing
\[
\psi_{x}(t): = \frac{1}{2}\left\{f(x+t)+f(x-t)-2f(x)\right\},
\]
we get
\[
T_{n,A}(f;x)-f(x) = \frac{2}{\pi}\int_{0}^{\pi}\psi_{x}(t)K_n(t)\,dt
\]
and
\begin{equation}\label{eq56}
\|T_{n,A}(f)-f\| = \mathcal{O}\left(\int_{0}^{a_{n,n}}|\psi_{x}(t)||K_n(t)|\,dt+
\int_{a_{n,n}}^{\pi}|\psi_{x}(t)||K_n(t)|\,dt\right).
 \end{equation}

By the well-known inequality $\sin\theta \geq \frac{2}{\pi}\theta $ for $ 0 \leq \theta \leq \frac{\pi}{2} $, $|\sin\theta |\leq 1$, and Lemma~\ref{theorem7} we get
\begin{align}\label{eq57}
\int_{0}^{a_{n,n}}|\psi_{x}(t)||K_n(t)|\,dt &= \mathcal{O}(1)\int_{0}^{a_{n,n}}t^{-1}\omega(t)dt\nonumber\\
&= \mathcal{O}\left(a_{n,n}H\left(a_{n,n}\right)\right).
\end{align}

Now using Lemma~\ref{theorem9} and~\eqref{eq13}, we have
\begin{align}\label{eq58}
\int_{a_{n,n}}^{\pi}|\psi_{x}(t)||K_n(t)|\,dt &= \mathcal{O}\left(a_{n,n}\right)\int_{a_{n,n}}^{\pi}t^{-2}\omega(t)dt\nonumber \\
&= \mathcal{O}\left(a_{n,n}H\left(a_{n,n}\right)\right).
\end{align}

Substituting~\eqref{eq57} and~\eqref{eq58} into~\eqref{eq56}, this implies~\eqref{eq55}.

The proof is now complete.
\end{proof}

\begin{theorem}\label{theorem11}
Let $\{a_{n,k}\}$ satisfy the conditions~\eqref{eq11} and
\[
\sum_{k = 0}^{m-1}{(k+1)}^{\beta}\left|\frac{a_{n,k}}{{(k+1)}^{\beta}}-\frac{a_{n,k+1}}{{(k+2)}^{\beta}}\right| = \mathcal{O}(a_{n,m}),
\]
for all $m = 0,1,\dots,n;\, (n = 0,1,\dots),$ and~\eqref{eq13}.

Then
\begin{equation}\label{eq59}
\|T_{n,A}(f)-f\|
 = \mathcal{O}\left(\omega\left(\frac{\pi}{n+1}\right)+a_{n,n}H\left(\frac{\pi}{n+1}\right)\right).
\end{equation}

If, in addition, $\omega (t)$ satisfies~\eqref{eq14} then
\begin{equation}\label{eq60}
\|T_{n,A}(f)-f\| = \mathcal{O}\left(a_{n,n}H\left(\frac{\pi}{n+1}\right)\right).
\end{equation}
\end{theorem}

\begin{proof}
Reasoning as in the proof of Theorem~\ref{theorem10} we have

\begin{align}\label{eq61}
\|T_{n,A}-f\| &= \mathcal{O}\left(\frac{2}{\pi}\right)\int_{0}^{\pi}|\psi_{x}(t)||K_n(t)|\,dt\nonumber \\
&= \mathcal{O}\left(\frac{2}{\pi}\right)\left(\!\int_{0}^{\frac{\pi}{n+1}}|\psi_{x}(t)||K_n(t)|\,dt+\!\int_{\frac{\pi}{n+1}}^{\pi}|\psi_{x}(t)||K_n(t)|\,dt\right).
\end{align}

Using the inequalities $\sin\theta \geq \frac{2}{\pi}\theta $ for $ 0 \leq \theta \leq \frac{\pi}{2} $, $\displaystyle |\sin t|\leq t$ and~\eqref{eq11} we have

\begin{align}\label{eq62}
\int_{0}^{\frac{\pi}{n+1}}|\psi_{x}(t)||K_n(t)|\,dt &= \mathcal{O}(n+1)\int_{0}^{\frac{\pi}{n+1}}\omega(t)\,dt\nonumber \\
&= \mathcal{O}\left(\omega\left(\frac{\pi}{n+1}\right)\right).
\end{align}

Lemma~\ref{theorem9} and~\eqref{eq13} implies
\begin{align}\label{eq63}
\int_{\frac{\pi}{n+1}}^{\pi}|\psi_{x}(t)||K_n(t)|\,dt &= \mathcal{O}\left(a_{n,n}\right) \int_{\frac{\pi}{n+1}}^{\pi}t^{-2}\omega(t)\,dt\nonumber \\
&= \mathcal{O}\left(a_{n,n}H\left(\frac{\pi}{n+1}\right)\right).
\end{align}

Putting~\eqref{eq62} and~\eqref{eq63} into~\eqref{eq61} we obtain~\eqref{eq59}.

Now we prove~\eqref{eq60}. It is clear that for $0\leq l\leq n-1$ and by our assumption we have
\begin{align*}
\mathcal{O}(a_{n,n}) &= \sum_{k = 0}^{n-1}{(k+1)}^{\beta}\left|\frac{a_{n,k}}{{(k+1)}^{\beta}}-\frac{a_{n,k+1}}{{(k+2)}^{\beta}}\right|\nonumber \\
&\geq {(l+1)}^{\beta}\sum_{k = l}^{n-1}\left|\frac{a_{n,k}}{{(k+1)}^{\beta}}-\frac{a_{n,k+1}}{{(k+2)}^{\beta}}\right|\\
&\geq {(l+1)}^{\beta}\left|\frac{a_{n,l}}{{(l+1)}^{\beta}}-\frac{a_{n,n}}{{(n+1)}^{\beta}}\right|\\
&\geq {(l+1)}^{\beta}\left|\frac{a_{n,l}}{{(l+1)}^{\beta}}-\frac{a_{n,n}}{{(l+1)}^{\beta}}\right|\geq a_{n,l}-a_{n,n},
\end{align*}
from which we conclude that
\[
a_{n,l} = \mathcal{O}\left(a_{n,n}\right).
\]

Using the last relation we get
\begin{align}\label{eq64}
1 = \sum_{l = 0}^{n}a_{n,l} &= \mathcal{O}\left(1\right)\sum_{l = 0}^{n}a_{n,n}
= \mathcal{O}\left((n+1)a_{n,n}\right)\nonumber \\
&\Longrightarrow \frac{1}{n+1} = \mathcal{O}\left(a_{n,n}\right).
\end{align}

By Lemma~\ref{theorem6},~\eqref{eq64} and~\eqref{eq62} we obtain
\begin{align}\label{eq65}
\int_{0}^{\frac{\pi}{n+1}}|\psi_{x}(t)||K_n(t)|\,dt &= \mathcal{O}(n+1)\int_{0}^{\frac{\pi}{n+1}}\omega(t)\,dt\nonumber \\
&= \mathcal{O}\left(n+1\right){(n+1)}^{-2}H\left(\frac{\pi}{n+1}\right)\nonumber \\
&= \mathcal{O}\left({(n+1)}^{-1}H\left(\frac{\pi}{n+1}\right)\right)\nonumber \\
&= \mathcal{O}\left(a_{n,n}H\left(\frac{\pi}{n+1}\right)\right).
\end{align}
Thus~\eqref{eq63} and~\eqref{eq65} imply~\eqref{eq60}, which completes the proof.
\end{proof}

\begin{theorem}\label{theorem12}
Let us assume that~\eqref{eq11} and
\[
\sum_{k = m}^{\infty}{(k+1)}^{\beta}\left|\frac{a_{n,k}}{{(k+1)}^{\beta}}-\frac{a_{n,k+1}}{{(k+2)}^{\beta}}\right| = \mathcal{O}(a_{n,m}),
\]
for all $m = 0,1,\dots,n;\, (n = 0,1,\dots)$.
Then
\begin{equation}\label{eq67}\nonumber
\|T_{n,A}(f)-f\| = \mathcal{O}\left(\omega\left(\frac{\pi}{n+1}\right)+\sum_{v = 1}^{n}v^{-1}\omega\left(\frac{\pi}{v}\right)\sum_{j = 0}^{v}a_{n,j}\right).
\end{equation}
\end{theorem}

\begin{proof}
Once again, we start from~\eqref{eq61}. For the first integral we use~\eqref{eq62}:
\begin{equation}\label{eq68}
\int_{0}^{\frac{\pi}{n+1}}|\psi_{x}(t)||K_n(t)|\,dt
= \mathcal{O}\left(\omega\left(\frac{\pi}{n+1}\right)\right).
\end{equation}

Applying Lemma~\ref{theorem9} we obtain
\begin{align*}\label{eq69}
\int_{\frac{\pi}{n+1}}^{\pi}|\psi_{x}(t)||K_n(t)|\,dt
&= \mathcal{O}\left(\sum_{v = 1}^{n}\int_{\frac{\pi}{v+1}}^{\frac{\pi}{v}}t^{-1}\omega(t)A_{n,\tau}dt\right)\nonumber \\
&= \mathcal{O}\left(\sum_{v = 1}^{n}(v+1)\omega\left(\frac{\pi}{v}\right)\sum_{j = 0}^{v+1}a_{n,j}\int_{\frac{\pi}{v+1}}^{\frac{\pi}{v}}dt\right)\nonumber \\
&= \mathcal{O}\left(\sum_{v = 1}^{n}v^{-1}\omega\left(\frac{\pi}{v}\right)\sum_{j = 0}^{v}a_{n,j}\right),
\end{align*}
since
\[
\int_{\frac{\pi}{v+1}}^{\frac{\pi}{v}}dt = \frac{\pi}{v(v+1)}.
\]
This with~\eqref{eq68} and~\eqref{eq61} completes the proof.
\end{proof}

\begin{theorem}\label{theorem13}
Let $\{a_{n,k}\}$ satisfy the conditions~\eqref{eq11} and
\[
\sum_{k = m}^{\infty}{(k+1)}^{\beta}\left|\frac{a_{n,k}}{{(k+1)}^{\beta}}-\frac{a_{n,k+1}}{{(k+2)}^{\beta}}\right| = \mathcal{O}(a_{n,m}),
\]
for all $m = 0,1,\dots,n;\, (n = 0,1,\dots),$
and suppose that $\omega (t)$ is such that~\eqref{eq13} and~\eqref{eq14} hold. Then
\begin{equation}\label{eq70}
\|T_{n,A}(f)-f\| = \mathcal{O}\left(a_{n,0}H\left(a_{n,0}\right)\right).
\end{equation}
\end{theorem}

\begin{proof}
We have seen in the proof of Lemma~\ref{theorem9} that
\begin{equation*}
\left|K_n(t)\right|
\leq \frac{\pi^2 a_{n,n}}{t^2}
+\frac{\pi^2}{t^2}\sum_{k = 0}^{n-1}{(k+1)}^{\beta}\left|\frac{a_{n,k}}{{(k+1)}^{\beta}}-\frac{a_{n,k+1}}{{(k+2)}^{\beta}}\right|.
\end{equation*}

Based on our assumption, for $n\geq 0$, we have
\begin{align*}
a_{n,n} &= \sum_{k = n}^{\infty}{(k+1)}^{\beta}\left|\frac{a_{n,k}}{{(k+1)}^{\beta}}-\frac{a_{n,k+1}}{{(k+2)}^{\beta}}\right|\\
&\leq \sum_{k = 0}^{\infty}{(k+1)}^{\beta}\left|\frac{a_{n,k}}{{(k+1)}^{\beta}}-\frac{a_{n,k+1}}{{(k+2)}^{\beta}}\right| = \mathcal{O}(a_{n,0}),
\end{align*}
and thus
\begin{equation*}
\left|K_n(t)\right|
\leq \frac{\pi^2 a_{n,0}}{t^2}
+\frac{\pi^2}{t^2}\sum_{k = 0}^{\infty}{(k+1)}^{\beta}\left|\frac{a_{n,k}}{{(k+1)}^{\beta}}-\frac{a_{n,k+1}}{{(k+2)}^{\beta}}\right|
= \mathcal{O}\left(\frac{a_{n,0}}{t^2}\right).
\end{equation*}

Now we write
\begin{equation}\label{eq72}
\|T_{n,A}(f)-f\| = \mathcal{O}\left(\int_{0}^{a_{n,0}}|\psi_{x}(t)||K_n(t)|\,dt+\int_{a_{n,0}}^{\pi}|\psi_{x}(t)||K_n(t)|\,dt\right).
\end{equation}

Similarly, as in~\eqref{eq57}), we have
\begin{align}\label{eq73}
\int_{0}^{a_{n,0}}|\psi_{x}(t)||K_n(t)|\,dt &= \mathcal{O}(1)\int_{0}^{a_{n,0}}t^{-1}\omega(t)\,dt\nonumber\\
&= \mathcal{O}\left(a_{n,0}H\left(a_{n,0}\right)\right).
\end{align}

Moreover, as in~\eqref{eq58}, we obtain
\begin{align}\label{eq74}
\int_{a_{n,0}}^{\pi}|\psi_{x}(t)||K_n(t)|\,dt &= \mathcal{O}\left(a_{n,0}\right)\int_{a_{n,0}}^{\pi}t^{-2}\omega(t)\,dt\nonumber \\
&= \mathcal{O}\left(a_{n,0}H\left(a_{n,0}\right)\right).
\end{align}

Putting~\eqref{eq73} and~\eqref{eq74} into~\eqref{eq72}, we get~\eqref{eq70}, which completes the proof.
\end{proof}

\begin{remark}
Based on the remark given at the end of the first section (of this paper), we conclude that our results hold true if we assume that either the conditions $\{{(k+1)}^{-\beta}a_{n,k}\}\in HBVS$ or $\{{(k+1)}^{-\beta}a_{n,k}\}\in RBVS$,
in the appropriate cases.
\end{remark}

\begin{remark}
If we put $\beta = 0$ in our theorems, then we exactly obtain all results proved in~\cite{L}.
\end{remark}

\section{Corollaries}

Assume that $A: = (a_{n,k})$ is a matrix defined by
\begin{equation*}
    a_{n,k} =
\begin{cases}
               \frac{p_{n-k}}{P_{n}}, & 0\leq k\leq n\\
               0, & \text{otherwise},
           \end{cases}
\end{equation*}
where $\{p_{k}\}_{k = 0,1,\dots,n}$, is a sequence of non-negative numbers with $P_{n}: = \sum_{k = 0}^{n}p_{k} \ne 0$. In this case, the transform $T_{n,A}(f;x)$ reduces to polynomials of the form
\[
T_{n,A}(f;x)\equiv N_n(f;x) = \frac{1}{P_{n}}\sum_{k = 0}^{n}p_{n-k}s_k(f;x).
\]

Consequently, we have the following.
\begin{corollary}\label{cor41}
Let all conditions of Theorem~\ref{theorem13} be satisfied. Then
\begin{equation}
\nonumber
\|N_{n}(f)-f\|
= \mathcal{O}\left(\frac{p_{n}}{P_{n}}H\left(\frac{p_{n}}{P_{n}}\right)\right).
\end{equation}
\end{corollary}

If the matrix $A: = (a_{n,k})$ is defined as follows:
\begin{equation*}
    a_{n,k} =
\begin{cases}
               \frac{p_{k}}{P_{n}}, & 0\leq k\leq n\\
               0, & \text{otherwise},
           \end{cases}
\end{equation*}
where $\{p_{k}\}_{k = 0,1,\dots,n}$, is a sequence of non-negative numbers with $P_{n}: = \sum_{k = 0}^{n}p_{k} \ne 0$. Now, the transform $T_{n,A}(f;x)$ are the polynomials
\[
T_{n,A}(f;x)\equiv R_n(f;x) = \frac{1}{P_{n}}\sum_{k = 0}^{n}p_{k}s_k(f;x).
\]

Therefore, we have:
\begin{corollary}\label{cor42}
Let all conditions of Theorem~\ref{theorem10} be satisfied. Then
\begin{equation*}
\|R_{n}(f)-f\|
= \mathcal{O}\left(\frac{p_{n}}{P_{n}}H\left(\frac{p_{n}}{P_{n}}\right)\right).
\end{equation*}
\end{corollary}

\begin{remark}
For $\beta = 0$, Corollary~\ref{cor41} and Corollary~\ref{cor42} reduce to two corollaries given in~\cite{L} (see pages 108-109).
\end{remark}

If $\omega(u) = u^{\alpha}$, $0<\alpha \leq 1$, and $f\in \text{Lip}(\omega)\equiv \text{Lip}(\alpha)$, then
$\omega(t) = \mathcal{O}\left(\delta^{\alpha}\right),\, \delta \geq 0$. Hence,
if $f\in \text{Lip}(\alpha)$ and
\begin{equation*}
   H(u) =
\begin{cases}
               u^{\alpha -1}, & 0< \alpha <1;\\
               \log \left(\frac{\pi}{u}\right), & \alpha = 1,
           \end{cases}
\end{equation*}
and then from Corollary~\ref{cor42} we get:

\begin{corollary}\label{cor43}
Let $f\in C[0,2\pi]$, $f\in \text{Lip}(\alpha)$, $(0<\alpha \leq 1)$, and let $\{p_k/P_n\}$ satisfies the condition $\eqref{star}$.
Then the degree of approximation of $f$ by $R_{n}(f;x)$
of its Fourier series is given by
\begin{equation*}
   \|R_{n}(f)-f\|
=
\begin{cases}
               \mathcal{O}\left({\left(\frac{p_{n}}{P_{n}}\right)}^{\alpha}\right), & 0< \alpha <1;\\
               \mathcal{O}\left(\frac{p_{n}}{P_{n}}\log \left(\frac{\pi P_{n}}{p_{n}}\right)\right), & \alpha = 1.
           \end{cases}
\end{equation*}
\end{corollary}

\begin{remark}
Note that for $\beta = 0$ in Corollary~\ref{cor43} we obtain one of the theorems proved in~\cite{CH4}.
\end{remark}

\subsection*{Acknowledgements}
The author expresses his great gratitude to the anonymous referee who with her/his report draw my attention that Lemma~\ref{theorem6} is a consequence of Lemma~\ref{theorem7} offering a simple proof, and also made several comments which definitely improved this paper. Also, I would like to thank the editor of Communications in Mathematics who helped me with many advises pertaining to the language used throughout the paper.

\EditInfo{%
    13 April 2020}{%
    10 July 2020}{%
    Karl Dilcher}

\end{paper}